\documentclass[12pt,leqno]{article}

\def\doublespace{\normalbaselineskip=24pt\baselineskip=\normalbaselineskip}

\newtheorem{theorem}{Theorem}
\newtheorem{lemma}{Lemma}
\newtheorem{corollary}{Corollary}

\begin{document}

\title{\textbf{A Convergence Analysis on URV Refinement}}
\author{\Large Limin Wu\\
               Department of Applied Mathematics\\
               Florida Institute of Technology\\
               Melbourne, FL. 32901-6988\\}
\date{August 1997}
\maketitle

\vspace{0.5in}

\begin{abstract}
Recently, Stewart gave an algorithm for computing a rank revealing
URV decomposition of a rectangular matrix. His method makes use
of a refinement iteration to achieve an improved estimate of
the smallest singular value and its corresponding singular vectors of the
matrix. Here, a new proof is given for the convergence of the 
refinement iteration. This analysis is carried out under slightly
weaker assumptions than those of Mathias and Stewart.
\end{abstract}


{\doublespace             

\textbf{1. Introduction.}

In \cite{Stew1}, Stewart gave an updating algorithm for subspace tracking.
His algorithm makes use of a refinement iteration, called URV refinement
in the literature, to achieve an improved
estimate of the smallest singular value and its corresponding singular
vectors of a nonsingular upper triangular matrix. The URV refinement 
can be briefly described as follows.

Consider a real $n\times n$ nonsingular upper triangular matrix $R$.
Let $R^{(0)} = R$  be partitioned as
\begin{equation}\label{eq1_1}
  R^{(0)} =
  \left[ \begin{array}{cc}
          S^{(0)} & h^{(0)}\\
          0       & e^{(0)}
          \end{array} \right],
\end{equation}
where $S^{(0)}$ is an $(n-1)\times (n-1)$ upper triangular matrix, 
$h^{(0)}$ is an $(n-1)$-vector, and $e^{(0)}$ is a scalar. 
Then a sequence of orthogonal
matrices, $Q^{(1)}, Q^{(2)}, \cdots, Q^{(2k-1)}, Q^{(2k)}$,
each determined as products of Givens rotations, 
is constructed such 
that, for $k \ge 1$, 
\begin{equation}\label{eq1_2}
  R^{(2k-1)} \equiv R^{(2k-2)}[Q^{(2k-1)}]^T =
  \left[ \begin{array}{cc}
          S^{(2k-1)} & 0\\
          h^{(2k-1)} & e^{(2k-1)}
          \end{array} \right],
\end{equation}
\begin{equation}\label{eq1_3}
  R^{(2k)} \equiv Q^{(2k)}R^{(2k-1)} =
  \left[ \begin{array}{cc}
          S^{(2k)} & h^{(2k)}\\
          0        & e^{(2k)}
          \end{array} \right],
\end{equation}
where $S^{(2k-1)}, S^{(2k)}$ are $(n-1)\times(n-1)$ upper triangular
matrices. 

The URV refinement is identified by Chandrasekaran and Ipsen 
\cite{ChIp} as an incomplete version of the QR algorithm for computing 
the singular value decomposition of an upper triangular matrix.
Stewart and Mathias \cite{Stew2,MaSt} discussed the URV refinement 
in a broader framework of block QR iterations, where 
$S^{(l)}$ are allowed to be $k\times k\, (1 \le k < n)$ matrices,
not necessarily upper triangular, and the
$e^{(l)}$ are then $(n-k) \times (n-k)$ matrices.
They established error bounds and derived
convergence properties for the singular values of $S^{(l)}$ and $e^{(l)}$.
In particular, for the special case considered in this paper, 
they proved that, if $|e^{(0)}|/\sigma_{min}(S^{(0)}) < 1$, 
then the URV refinement computes the smallest singular value. 
We have used $\sigma_{min}(\cdot)$
to denote the smallest singular value of a matrix.
We will also use $\sigma_i(\cdot)$ to denote the $ith$ largest singular
value of a matrix and $\|\cdot\|$ to denote the 2-norm of a matrix 
throughout the paper.

To facilitate comparison with the new convergence proof given here, 
we restate a theorem from \cite{MaSt} for the case $k = n-1$.
\begin{theorem}[Mathias and Stewart, 1993]\label{thm1}
Let $S^{(l)}, e^{(l)}$, and $h^{(l)}$ be defined as in (1)-(3).
For $l \ge 1$, we have
\begin{enumerate}
  \item $|e^{(l)}| \le |e^{(l-1)}|$;
  \item $\sigma_j(S^{(l)}) \ge \sigma_j(S^{(l-1)}),\;j = 1,\ldots,n-1$;
  \item $\|h^{(l)}\| \le \rho^{(l)}\cdots\rho^{(0)}\|h^{(0)}\| \le 
        (\rho^{(0)})^l\|h^{(0)}\|$, 
        where $\rho^{(l)} \equiv |e^{(l)}|/\sigma_{min}(S^{(l)})$;

and if $\rho^{(0)} < 1$, then
  \item $\lim_{l\rightarrow\infty}|e^{(l)}| = \sigma_n(R)$;
  \item $\lim_{l\rightarrow\infty}\sigma_j(S^{(l)}) = \sigma_j(R),\; 
        j = 1,\ldots,n-1$.
\end{enumerate}
\end{theorem}

The assumption $\rho^{(0)} < 1$ is needed for the method of proof used to 
establish parts $4-5$ of the above theorem, but is not a necessary condition
for the convergence of the algorithm. An example which illustrates this
fact is
\[ R = \left[ \begin{array}{ccl}
              1 & 0 & 10^{-6} \\
              0 & 2 & 10^{-6} \\
              0 & 0 & 10
       \end{array} \right].
\]
For this example, a MATLAB implementation of the URV refinement 
 yields an approximation to the smallest singular value of $R$ after
$14$ iterations as $e^{(28)} = 9.9\cdots948e\!-\!01$ in double precision. 
This is very close to the
smallest singular value of $R$, $s_3 = 9.9\cdots950e\!-\!01$ 
computed using the MATLAB SVD routine.

\textbf{2. Convergence Analysis.}

The singular value decomposition (SVD) of $R^{(l)}$'s provides a basis for
our convergence analysis.
Let $R^{(0)}$ have the SVD 
$R^{(0)} = U^{(0)} \Sigma [V^{(0)}]^T$, where 
$\Sigma = diag(\sigma_1,\sigma_2, \cdots,\sigma_n)$, with 
\[ \sigma_1~\geq~\sigma_2~\geq~\cdots~\geq~\sigma_n>0.\] 
Defining $V^{(1)} = Q^{(1)}V^{(0)}$ and $U^{(2)} = Q^{(2)}U^{(0)}$, 
then $R^{(1)}$ and $R^{(2)}$ have  SVD's,
$R^{(1)} = U^{(0)} \Sigma [V^{(1)}]^T$ and
$R^{(2)} = U^{(2)} \Sigma [V^{(1)}]^T$, respectively.
For $k \geq 2$, define 
\begin{equation}\label{eq21}
  V^{(2k-1)} = Q^{(2k-1)}V^{(2k-3)} = G^{(2k-1)}V,\; \textrm{where} \;
  G^{(2k-1)} = Q^{(2k-1)}\cdots Q^{(1)},
\end{equation}
and
\begin{equation}\label{eq51}
  U^{(2k)} = Q^{(2k)}U^{(2k-2)} = G^{(2k)}U,\; \textrm{where}\;
  G^{(2k)} = Q^{(2k)}\cdots Q^{(2)}.
\end{equation}
Then $R^{(2k-1)}$ and $R^{(2k)}$ have SVD's:
\begin{equation}\label{eq6}
  R^{(2k-1)} =  U^{(2k-2)} \Sigma [V^{(2k-1)}]^T,
\end{equation}
and 
\begin{equation}\label{eq7}
  R^{(2k)} = U^{(2k)} \Sigma [V^{(2k-1)}]^T.
\end{equation}  

Denote $R^{(l)} = [r_{ij}^{(l)}], V^{(2k-1)} = [v_{ij}^{(2k-1)}]$, and
$U^{(2k)} = [u_{ij}^{(2k)}]$; then $r^{(l)}_{nn} = e^{(l)}$.
Let $v^{(0)}$ and $u^{(0)}$ be the last columns of matrices $V^{(0)}$ 
and $U^{(0)}$, 
respectively. Let $g^{(2k-1)}$ and $g^{(2k)}$ contain the last rows of
$G^{(2k-1)}$ and $G^{(2k)}$, respectively. 
Then the following theorem, first given in \cite{Wu}, holds.
\begin{theorem}\label{thm2}
  Assuming in the URV refinement that $r^{(l)}_{nn}$ are kept positive,
then we have:
\begin{enumerate}
   \item if $v^{(0)}_{nn} \ne 0$, then $r^{(l)}_{nn}$ 
         converges to $\sigma_n$ monotonically; 
   \item if $v^{(0)}_{nn} \ne 0$ and $\sigma_{n-1} > \sigma_n$, then
         $\langle g^{(2k-1)},v^{(0)}\rangle \equiv 
         [g^{(2k-1)}]^Tv^{(0)}$ and
         $\langle g^{(2k)},u^{(0)}\rangle \equiv 
         [g^{(2k)}]^Tu^{(0)}$ 
         converge to $\pm1$ monotonically;
   \item if $\langle g^{(2k-1)},v^{(0)}\rangle \rightarrow 1$, 
         then $v^{(0)}_{nn} \ne 0$.
\end{enumerate}
\end{theorem}
The condition $\sigma_{n-1} > \sigma_n$ in part 2 of the theorem says
that the smallest singular value $\sigma_n$ is not repeated, that is, it is
simple. Therefore it 
has unique left and right singular vectors associated with it.
To prove this theorem, we need the following lemmas.
\begin{lemma}\label{lem1}
The smallest singular value of a square nonsingular triangular matrix is 
not greater than the absolute value of any diagonal element of the matrix.
\end{lemma}
\textbf{Proof.} see Lawson and Hanson \cite[p.29, (6.3)]{LaHa}. 
\rule{2mm}{3mm}

\begin{lemma}\label{lem2}
The sequence $\{|r_{nn}^{(l)}|\}_1^\infty$ obtained from the URV
refinement is nonincreasing and converges. In particular, 
$\{r_{nn}^{(l)}\}_1^\infty$ is nonincreasing and converges
if $r_{nn}^{(l)}$ are kept positive in the URV refinement.
\end{lemma}
\textbf{Proof.} Since orthogonal matrices preserve the 2-norm of vectors, 
we have
\[ \|\left[ \begin{array}{c}
	 r^{(0)}_{1n} \\ \vdots \\ r^{(0)}_{nn}
     \end{array} \right]\| \geq r^{(0)}_{nn} = 
     \|(r_{n1}^{(1)},\ldots,r_{nn}^{(1)})^T\|_2 \geq |r_{nn}^{(1)}| =
   \|\left[ \begin{array}{c} 
     r_{1n}^{(2)} \\ \vdots \\ r_{nn}^{(2)}
   \end{array} \right]\|_2 \geq |r_{nn}^{(2)}| = \cdots > 0.
\]
Thus, $\{|r_{nn}^{(l)}|\}_1^\infty$ is nonincreasing and bounded
below by $0$. It follows that this sequence has a limit. 
If we choose Givens rotations in the refinement process in 
such a way that $r_{nn}^{(l)}$ are kept positive, then the sequence
$\{r_{nn}^{(l)}\}_1^\infty$ has a limit. \rule{2mm}{3mm}

\begin{lemma}\label{lem3}
If $r_{nn}^{(l)}$ are kept positive in the URV refinement, then
\begin{enumerate}
  \item $u^{(0)}_{nn} = v_{nn}^{(2k-1)} = u_{nn}^{(2k)} = 0$, 
        for $k \ge 1$, provided $v^{(0)}_{nn} = 0$.
  \item $u^{(0)}_{nn} > 0$, $v_{nn}^{(2k-1)} > 0$, and $u_{nn}^{(2k)} >
        0$, for $k \geq 1$, provided $v^{(0)}_{nn} > 0$.
  \item $u^{(0)}_{nn} < 0$, $v_{nn}^{(2k-1)} < 0$, and $u_{nn}^{(2k)} <
	0$, for $k \geq 1$, provided $v^{(0)}_{nn} < 0$.
\end{enumerate}
\end{lemma}
\textbf{Proof.} 
Write the SVD of $R^{(0)}$ as 
\begin{equation}\label{eq71}
  R^{(0)}V^{(0)} = U^{(0)}\Sigma.
\end{equation}
Equating the corner elements at the $(n,n)$-position on both
sides (\ref{eq71}) gives 
\begin{equation}\label{eq72}
  r^{(0)}_{nn}v^{(0)}_{nn} = \sigma_nu^{(0)}_{nn}.
\end{equation}
Also, SVD (\ref{eq6}) can be written as
\begin{equation}\label{eq8}
  [U^{(2k-2)}]^T R^{(2k-1)} = \Sigma[V^{(2k-1)}]^T.
\end{equation}
Since $R^{(2k-1)}$ is of form (\ref{eq1_2}), it is easy to see that
\begin{equation}\label{eq9}
  u_{nn}^{(2k-2)}r_{nn}^{(2k-1)} = \sigma_n v_{nn}^{(2k-1)},\;
  \textrm{for}\; k \geq 1.
\end{equation}
Similarly, writing (\ref{eq7}) as
\begin{equation}\label{eq10}
  R^{(2k)} V^{(2k-1)} = U^{(2k)} \Sigma,
\end{equation}
we have
\begin{equation}\label{eq11}
  r_{nn}^{(2k)}v_{nn}^{(2k-1)} = u_{nn}^{(2k)} \sigma_n, \;
  \textrm{for}\; k \geq 1.
\end{equation}
Since we have assumed that $r_{nn}^{(l)} > 0$ and $\sigma_n > 0$, the
conclusions are easily drawn using equations (\ref{eq72}),
(\ref{eq9}), and (\ref{eq11}). \rule{2mm}{3mm}

\begin{lemma}\label{lem4}
If $v^{(0)}_{nn} \ne 0$ and  $r_{nn}^{(l)}$ are kept positive
in the URV refinement, then $\{v_{nn}^{(2k-1)}\}_1^\infty$ and
$\{u_{nn}^{(2k)}\}_1^\infty$ converge monotonically to the same nonzero 
limit.
\end{lemma}
\textbf{Proof.}
We first assume $v^{(0)}_{nn} > 0$.  According to 
Lemma~\ref{lem3} , 
$u_{nn}^{(2k)}$ are also positive for $k \geq 1$. By manipulating 
(\ref{eq9}) and (\ref{eq11}) we obtain
\begin{equation}\label{eq12}
  \frac{u_{nn}^{(2k-2)}}{u_{nn}^{(2k)}} =
  \frac{\sigma_n^2}{r_{nn}^{(2k-1)}r_{nn}^{(2k)}}.
\end{equation}
By Lemma~\ref{lem1} and Lemma~\ref{lem2}, 
we have $r_{nn}^{(2k)} \geq \sigma_n$ and 
$r_{nn}^{(2k-1)} \geq r_{nn}^{(2k)}$. It follows that
the right hand side of (\ref{eq12}) is less than or 
equal to one. Therefore,
$\{u_{nn}^{(2k)}\}_1^\infty$ is a nondecreasing sequence. The
orthogonality of $U^{(2k)}$ means that $u_{nn}^{(2k)}$ is bounded from 
above by
one. Hence, $\{u_{nn}^{(2k)}\}_1^\infty$ has a positive limit. 
Also, the relation (\ref{eq9}) tells us that
$\{v_{nn}^{(2k-1)}\}_1^\infty$ has the same limit as
$\{u_{nn}^{(2k)}\}_1^\infty$ does.
For the case $v^{(0)}_{nn} < 0$ the proof is similar. \rule{2mm}{3mm}

\noindent \textbf{Proof of Theorem 2.}

\noindent Part $1$.
Since both $r_{nn}^{(l)}$ and $u_{nn}^{(2k)}$ converge,
taking the limit
on both sides of (\ref{eq12}) yields
$\lim_{l\rightarrow\infty}r_{nn}^{(l)} = \sigma_n$ and the convergence is
monotone by Lemma 2.

\noindent Part $2$. 
By (\ref{eq21}) and (\ref{eq51}), 
$\langle g^{(2k-1)}, v^{(0)}\rangle = v_{nn}^{(2k-1)}$ and
$\langle g^{(2k)}, u^{(0)}\rangle = u_{nn}^{(2k)}$.
We prove that   
$\lim_{k\rightarrow\infty}v_{nn}^{(2k-1)} =
 \lim_{k\rightarrow\infty}u_{nn}^{(2k)} = \pm1$ under the assumption.
Suppose $\lim_{k\rightarrow\infty}u_{nn}^{(2k)} = a$. 
Apparently $|a| \le 1$. By Lemma 4, 
$\lim_{k\rightarrow\infty}v_{nn}^{(2k-1)} = a$. 
Equating the last rows in both sides of (\ref{eq10}) gives 
\[ r_{nn}^{(2k)}(v_{n1}^{(2k-1)},\ldots,v_{nn}^{(2k-1)}) =
   (\sigma_1 u_{n1}^{(2k)},\ldots,\sigma_n u_{nn}^{(2k)}).
\]
Taking the 2-norm of the above equation and then squaring both sides
gives
\[ (r_{nn}^{(2k)})^2 = 
   \sigma_1^2(u_{n1}^{(2k)})^2+\cdots +\sigma_n^2(u_{nn}^{(2k)})^2.
\]
Rewriting the above equation and considering the ordering of $\sigma_i$'s
we get
\begin{eqnarray*}
  &   & (r_{nn}^{(2k)})^2-\sigma_n^2 (u_{nn}^{(2k)})^2 \\
  & = & \sigma_1^2 (u_{n1}^{(2k)})^2+\cdots +\sigma_{n-1}^2 
 (u_{n,n-1}^{(2k)})^2 \\ 
 &\ge & \sigma_{n-1}^2 (u_{n1}^{(2k)})^2+\cdots+ 
  \sigma_{n-1}^2(u_{n,n-1}^{(2k)})^2 \\
  & = & \sigma_{n-1}^2 (1-(u_{nn}^{(2k)})^2) 
\end{eqnarray*}
Taking the limit on both sides of the above equation yields
\[ \sigma_n^2 (1-a^2) \ge \sigma_{n-1}^2 (1-a^2) \]
Since we have assumed $\sigma_{n-1} > \sigma_n$, the only way that this
inequality can hold is if $a = \pm1$. 

\noindent Part $3$. Since $v_{nn}^{(2k-1)} =  
\langle g^{(2k)}, v^{(0)}\rangle \rightarrow 1$,
in view of Lemma 3, it is obvious that $v^{(0)}_{nn} \ne 0$. 
\rule{2mm}{3mm}

\noindent \textbf{Note} It is a consequence of the 
standard theory of inner product space that 
\[\lim_{k\rightarrow\infty}\langle g^{(2k-1)}, v^{(0)}\rangle = 1\;\; 
\textrm{if and only if}\;\;
\lim_{k\rightarrow\infty}\|g^{(2k-1)}-v^{(0)}\| = 0.\]

Since $v^{(0)}_{nn} \ne 0$ is 
vital for the convergence of the URV refinement
when $\sigma_n$ is simple, 
it is desirable to know under what conditions the nonsingular
upper triangular matrix
$R$ has a simple smallest singular value and nonzero $v^{(0)}_{nn}$ 
in its SVD. 
A sufficient condition is given by the following theorem. In the proof of
the theorem
we drop the superscript $^{(0)}$ for $R, V$, and 
$U$ and related quantities.
Let $R_1$ be the matrix consisting of the first $n-1$ columns of $R$.
\begin{lemma}
  $\sigma_{n-1} \ge \sigma_{min}(R_1) \ge \sigma_n$.
\end{lemma}
\textbf{Proof.} see Lawson and Hanson \cite[p.26, (5.12)]{LaHa}.
\rule{2mm}{3mm}
\begin{theorem}
  If $\sigma_{min}(S) > \sigma_n$, then $\sigma_n$ is simple 
  and $v_{nn} \ne 0$.
\end{theorem}
\textbf{Proof.} 
Since $\sigma_{min}(S) = \sigma_{min}(R_1)$, 
it follows by Lemma 5 that $\sigma_n$ is simple.
To prove the second part, we will show that $v_{nn} = 0$ implies 
$\sigma_{min}(S) = \sigma_n$.
First, since $\sigma_{min}(S) = \sigma_{min}(R_1)$, the inequality
$\sigma_{min}(S) \ge \sigma_n$ follows from Lemma 5.
We now establish the reverse inequality.
Let $R$ have the SVD $R = U\Sigma V^T$. 
Let $u$ be the last column of $U$, $r$ the last column of $R$.
Write the SVD of $R$ as
\begin{equation}\label{eq13}
  U^TR = \Sigma V^T.
\end{equation}
Equating the corner elements at the (n,n) position on both sides of 
(\ref{eq13})
gives $u^Tr = \sigma_nv_{nn}$. The assumption of $v_{nn} = 0$ implies
$u^Tr = 0$. Now consider the equation
\begin{equation}\label{eq14}
  u^T RR^Tu = u^TU\Sigma^2U^Tu,
\end{equation}
or the equivalent form
\begin{equation}\label{eq15}
  (u^TR)(u^TR)^T = (u^TU)\Sigma^2(u^TU)^T = e_n^T\Sigma^2e_n,
\end{equation}
where $e_n$ is the unit vector with one in the last component.
Since $u^Tr = 0$, (\ref{eq15}) becomes 
\begin{equation}\label{eq16}
  (u^TR_1)(u^TR_1)^T = \sigma_n^2. 
\end{equation}
Letting $w = (u_{1n},\ldots,u_{n-1,n})^T$, we have $u^TR_1 = w^TS$. 
Thus (\ref{eq16}) can be further reduced to
\begin{equation}\label{eq17}
  (w^TS)(w^TS)^T = \sigma_n^2, \; \textrm{or}\;
  \|S^Tw\| = \sigma_n.
\end{equation}
By Lemma 3, $u_{nn} = v_{nn} = 0$, thus $\|w\| = \|u\| = 1$.
Therefore we have
\begin{equation}\label{eq18}
 \sigma_{min}(S) = \sigma_{min}(S^T) = \min_{\|x\| = 1}\|S^Tx\|
   \le \|S^Tw\| = \sigma_n.
\end{equation}
This completes the proof that $\sigma_{min}(S) = \sigma_n$. \rule{2mm}{3mm}
\begin{corollary}
  If $\sigma_{min}(S) > \sigma_{min}(R)$, then the 
URV refinement converges.
\end{corollary}

\noindent \textbf{Remark}
The assumption of $\sigma_{min}(S^{(0)}) > \sigma_{min}(R^{(0)})$ in 
the above corollary is weaker than the assumption 
$\sigma_{min}(S^{(0)}) > |e^{(0)}|$
used in Theorem 1 because $|e^{(0)}| \ge \sigma_{min}(R^{(0)})$. 
We may note, however, that if 
$\sigma_{min}(S^{(0)}) > \sigma_{min}(R^{(0)})$, then
the URV refinement will produce, for 
sufficiently large $l$, an $|e^{(l)}|$ and $\sigma_{min}(S^{(l)})$
such that $\sigma_{min}(S^{(l)}) > |e^{(l)}|$. This follows from
$\sigma_{min}(S^{(l)}) \ge \sigma_{min}(S^{(0)})$ (Theorem 1, Part 2) and
the fact that $\sigma_{min}(S^{(0)}) > \sigma_{min}(R^{(0)})$ implies
$|e^{(l)}| \rightarrow \sigma_{min}(R^{(0)})$ (Corollary 1).

\textbf{ }

\textbf{ }

\textbf{3. Conclusions}

In this paper we have shown that $v_{nn}^{(0)} \ne 0$ is sufficient for
the convergence of the sequence $\{ r_{nn}^{(l)} \}$ to the smallest
singular value of $R^{(0)}$ in the URV refinement (Theorem2, Part1). 
The following matrix 
\[ R^{(0)} = \left[ \begin{array}{ccl}
              1 & 0 & 0 \\
              0 & 9 & 1 \\
              0 & 1 & 10
       \end{array} \right]
\]
serves as a convenient example for which $v_{33}^{(0)} = 0$ and the sequence 
$\{ r_{33}^{(l)} \}$ fails to converge to $R^{(0)}$.
It is unknown whether $v_{nn}^{(0)} \ne 0$ is a necessary condition for 
$r_{nn}^{(l)} \rightarrow \sigma_{min}(R^{(0)})$. 

\textbf{Acknowledgments.}
The author thanks Professor Charles Fulton for helpful discussions 
and Professor Gary Howell for a suggestion leading to a simplified proof of
\linebreak Theorem \ref{thm2}.

} 

\end{document}